\documentclass{amsart}

\newtheorem{theorem}{Theorem}[section]

\newtheorem{proposition}[theorem]{Proposition}

\theoremstyle{definition}
\newtheorem{definition}[theorem]{Definition}

\theoremstyle{remark}
\newtheorem{remark}[theorem]{Remark}

\numberwithin{equation}{section}

\begin{document}


\title[Complete Vertical Graphs with CMC in Warped Products]{Complete Vertical Graphs with Constant Mean Curvature in Semi-Riemannian Warped Products}

\author{A. Caminha}
\address{Departamento de Matem\'atica, Universidade Federal do Cear\'a, Fortaleza,
Cear\'a, Brazil. 60455-760}
\email{antonio.caminha@gmail.com}

\author{H. F. de Lima}
\address{Departamento de Matem\'atica e Estat\'{\i}stica,
Universidade Federal de Campina Grande, Campina Grande, Para\'{\i}ba, Brazil. 58109-970}
\email{henrique@dme.ufcg.edu.br}

\thanks{We would like to thank professors L. Al\'{\i}as and A. G. Colares for having showed us their
preprint~\cite{Alias2:06}, which enabled us to set proposition~\ref{prop:Laplacian of h}.}

\thanks{The second author is partially supported by CAPES}


\date{September 20, 2006}


\begin{abstract}
In this paper we study complete vertical graphs of constant mean curvature
in the Hyperbolic and Steady State spaces. We first derive suitable
formulas for the Laplacians of the height function and of a support-like function
naturally attached to the graph; then, under appropriate restrictions on the values
of the mean curvature and the growth of the height function, we obtain necessary
conditions for the existence of such a graph. In the two-dimensional case we
apply this analytical framework to state and prove Bernstein-type results in each of
these ambient spaces.
\end{abstract}

\maketitle

\section{Introduction}

This paper deals with complete non-compact constant mean curvature graphs over a horosphere of the Hyperbolic space,
as well as over horizontal hyperplanes ({\em slices}) in the Steady State space. In connection with our work, L. Al\'{\i}as and
M. Dajczer (cf.~\cite{Alias:06}) studied properly immersed complete surfaces of the $3-$dimensional Hyperbolic space
contained between two horospheres, obtaining a Bernstein-type result for the case of constant mean curvature
between $-1$ and $1$. In de Sitter space, K. Akutagawa (cf.~\cite{Akutagawa:87}) proved that
complete spacelike hypersurfaces having constant mean curvature in a specific interval of the real line are
totally umbilical. Also for de Sitter space, among other interesting results S. Montiel (cf.~\cite{Montiel:03})
proves that, under an appropriate restriction on their Hyperbolic Gauss map, complete spacelike Hypersurfaces of
constant mean curvature greater than or equal to $1$ must actually have mean curvature $1$.

For the Lorentz case, our motivation to restrict attention to the Steady State space comes from the fact that there
exists a natural duality between the Gauss maps of Riemannian hypersurfaces of this space and those of the Hyperbolic space,
provided we model these as hyperquadrics of the Lorentz-Minkowski space
(cf. section~\ref{section:Vertical graphs in the Hyperbolic space}).
Besides, in physical context the Steady State space appears naturally as an exact solution for the Einstein equations,
being a cosmological model where matter is supposed to travel along geodesics normal to horizontal hyperplanes;
these, in turn, serve as the initial data for the Cauchy problem associated to those equations
(cf.~\cite{Hawking:73}, chapter $5$).

In this work we model both our ambient spaces as semi-Riemannian warped products to obtain necessary
conditions for the existence of the types of graphs mentioned in the beginning of this introduction.
More precisely, under appropriate restrictions on the values of the mean curvature and the growth of the height
function of these graphs, we actually prove that the mean curvature has to be identically $1$
(cf. Theorem~\ref{thm:de Sitter application} and Theorem~\ref{thm:Hyperbolic application}). We also prove
(under a slightly stronger hypothesis in the Hyperbolic case) that the scalar curvature of our graphs
cannot be globally bounded away from zero in a certain sense.
The analytical framework we use to prove the above-mentioned results consists of the generalized maximum
principle of Omori and Yau. Specifically, we apply lemma $3$ of~\cite{Akutagawa:87} on nonnegative solutions to
the partial differential inequality $\Delta g\geq ag^2$ ($a$ being a positive real constant) to a carefully chosen
combination of functions naturally attached to our immersions.

In dimension $2$, for complete surfaces of nonnegative Gaussian curvature, we are able to obtain Bernstein-type
theorems related to our previous general results by using the fact that those surfaces are parabolic in the sense of Riemann surfaces
({cf.~\cite{Huber:57}). Indeed, if the size of the gradient of the height function of the graph is suitably bounded,
then the graph has to be a horosphere in the $3-$dimensional Hyperbolic space
({cf. Theorem~\ref{thm:Hyperbolic Bernstein-type theorem}), or a horizontal plane in the
$3-$dimensional Steady State space ({cf. Theorem~\ref{thm:de Sitter Bernstein-type theorem}).

This paper is organized in the following manner: in section~\ref{section:Conformal vector fields} we discuss
general semi-Riemannian manifolds furnished with conformal vector fields, and derive a formula for the Laplacian
of a support-like function associated to an oriented Riemannian hypersurface of such an ambient space.
Section~\ref{section:Semi-Riemannian warped products} recasts the result of the previous one in the particular context
of semi-Riemannian warped products with Riemannian fiber; we also compute the Laplacian of a general height function
and close the section by defining the objects of our main interest, namely, vertical graphs over fibers of such an
warped product. Finally, sections~\ref{section:Vertical graphs in the steady state space} and
~\ref{section:Vertical graphs in the Hyperbolic space} are respectively devoted to applications of this
general picture to the special cases of the Steady State space and the Hyperbolic space.

Finally, it was communicated to us by professor L.J. Al\'{\i}as that our results about
spacelike surfaces in the 3-dimensional steady state space (Section 4)
are somewhat related to a work in progress due to him and A. L. Albujer
(cf.~\cite{Albujer:06}).

\section{Conformal vector fields}\label{section:Conformal vector fields}

Let $\overline M^{n+1}$ be a connected semi-Riemannian manifold with metric $\overline g=\langle\,\,,\,\,\rangle$
of index $\nu\leq 1$, and semi-Riemannian connection $\overline\nabla$. For a vector field $X\in\mathcal X(\overline M)$,
let $\epsilon(X)=\langle X,X\rangle$; $X$ is said to be a {\em unit} vector field if
$\epsilon(X)=\pm 1$, {\em timelike} if $\epsilon(X)=-1$.

A vector field $V$ on $\overline M^{n+1}$ is said to be {\em conformal} if
\begin{equation}
\mathcal L_V\langle\,\,,\,\,\rangle=2\phi\langle\,\,,\,\,\rangle
\end{equation}
for some function $\phi\in C^{\infty}(\overline M)$, where $\mathcal L$ stands for the Lie
derivative of the metric of $\overline M$. The function $\phi$ is called the
{\em conformal factor} of $V$.

Since $\mathcal L_V(X)=[V,X]$ for all $X\in\mathcal X(\overline M)$, it follows from the
tensorial character of $\mathcal L_V$ that $V\in\mathcal X(\overline M)$ is conformal if
and only if
\begin{equation}\label{eq:1.1}
\langle\overline\nabla_XV,Y \rangle+\langle X,\overline\nabla_YV\rangle=2\phi\langle X,Y\rangle,
\end{equation}
for all $X,Y\in\mathcal X(\overline M)$. In particular, $V$ is a Killing vector field
relatively to $\overline g$ if and only if $\phi\equiv 0$.

In all that follows, we consider {\em Riemannian immersions} $\psi:\Sigma^n\rightarrow\overline M^{n+1}$, namely, immersions
from a connected, $n-$dimensional orientable differentiable manifold $\Sigma$ into $\overline M$, such that the
induced metric $g=\psi^*(\overline g)$ turns $\Sigma$ into a Riemannian manifold
(in the Lorentz case $\nu=1$, we refer to $(\Sigma,g)$ as a {\em spacelike} hypersurface of $\overline M$),
with Levi-Civita connection $\nabla$. We orient $\Sigma$ by the choice of a unit normal vector field $N$ on it,
let $A$ denote the corresponding shape operator and $H=\epsilon(N)\,{\rm tr}(A)/n$ the corresponding mean curvature.

The following proposition appeared for the first time in~\cite{Sousa:04}, there in the Riemannian setting.
In a joint work with A. B. Barros and A. Brasil (cf.~\cite{Barros:05}) the first author generalized it to the Lorentz setting.
Here we present a unified version of it, together with a proof.

\begin{proposition}\label{prop:Laplacian of conformal vector field}
Let $\overline M^{n+1}$ be semi-Riemannian manifold furnished with a conformal vector
field $V$ with conformal factor $\phi:\overline M^{n+1}\rightarrow\mathbb R$, and
$\psi:\Sigma^n\rightarrow\overline M^{n+1}$ a Riemannian immersion. If $\eta=\langle V,N\rangle$, then
\begin{equation}\label{eq:Laplacian formula_I}
\Delta\eta=-\epsilon\,n\langle V,\nabla H\rangle-\epsilon\,\eta\left\{\overline{\rm{Ric}}(N,N)+|A|^2\right\}-n\left\{\epsilon\,H\phi+N(\phi)\right\},
\end{equation}
where $\epsilon=\epsilon(N)$, $\nabla H$ the gradient of $H$ in the metric of $\Sigma$,
$\overline{\rm{Ric}}$ is the Ricci tensor of $\overline M$ and $|A|$ is the Hilbert-Schmidt norm of $A$.
\end{proposition}

\begin{proof}
Fix $p\in\Sigma$ and let $\{e_k\}$ be an orthonormal moving frame on a neighborhood of $p$ in $\Sigma$,
geodesic at $p$. Extend the $e_k$ to a neighborhood of $p$ in $\overline M$, so that
$(\overline\nabla_Ne_k)(p)=0$, and let
$$V=\sum_l^n\alpha_le_l+\epsilon\,\eta N.$$
Then
\begin{eqnarray*}
\eta=\langle N,V\rangle\Rightarrow e_k(\eta)&=&\langle\overline\nabla_{e_k}N,V\rangle+\langle N,\overline\nabla_{e_k}V\rangle\\
&=&-\langle Ae_k,V\rangle+\langle N,\overline\nabla_{e_k}V\rangle,
\end{eqnarray*}
so that
\begin{eqnarray}\label{eq:I}
\Delta\eta&=&\sum_ke_k(e_k(\eta))=-\sum_ke_k\langle Ae_k,V\rangle+\sum_ke_k\langle N,\overline\nabla_{e_k}V\rangle\nonumber\\
&=&-\sum_k\langle\overline\nabla_{e_k}Ae_k,V\rangle-2\sum_k\langle Ae_k,\overline\nabla_{e_k}V\rangle+\sum_k\langle N,\overline\nabla_{e_k}\overline\nabla_{e_k}V\rangle.
\end{eqnarray}

Now, differentiating $Ae_k=\sum_lh_{kl}e_l$ with respect to $e_k$, one gets at $p$
\begin{eqnarray}\label{eq:II}
\sum_k\langle\overline\nabla_{e_k}Ae_k,V\rangle&=&\sum_{k,l}e_k(h_{kl})\langle e_l,V\rangle+\sum_{k,l}h_{kl}\langle\overline\nabla_{e_k}e_l,V\rangle\nonumber\\
&=&\sum_{k,l}\alpha_le_k(h_{kl})+\epsilon\sum_{k,l}h_{kl}\langle\overline\nabla_{e_k}e_l,N\rangle\langle V,N\rangle\nonumber\\
&=&\sum_{k,l}\alpha_le_k(h_{kl})+\epsilon\sum_{k,l}h_{kl}^2\eta\nonumber\\
&=&\sum_{k,l}\alpha_le_k(h_{kl})+\epsilon\,\eta|A|^2.
\end{eqnarray}

Asking further that $Ae_k=\lambda_ke_k$ at $p$ (which is always possible), we have at $p$
\begin{equation}\label{eq:III}
\sum_k\langle Ae_k,\overline\nabla_{e_k}V\rangle=\sum_k\lambda_k\langle e_k,\overline\nabla_{e_k}V\rangle=\sum_k\lambda_k\phi=\epsilon\,nH\phi.
\end{equation}

In order to compute the last summand of (\ref{eq:I}), note that the conformality of $V$ gives
$$\langle\overline\nabla_NV,e_k\rangle+\langle N,\overline\nabla_{e_k}V\rangle=0$$
for all $k$. Hence, differentiating the above relation in the direction of $e_k$, we get
$$\langle\overline\nabla_{e_k}\overline\nabla_NV,e_k\rangle+\langle\overline\nabla_NV,\overline\nabla_{e_k}e_k\rangle+\langle\overline\nabla_{e_k}N,\overline\nabla_{e_k}V\rangle+\langle N,\overline\nabla_{e_k}\overline\nabla_{e_k}V\rangle=0.$$
However, at $p$ one has
\begin{eqnarray*}
\langle\overline\nabla_NV,\overline\nabla_{e_k}e_k\rangle&=&\epsilon\langle\overline\nabla_NV,\langle\overline\nabla_{e_k}e_k,N\rangle N\rangle=\epsilon\langle\overline\nabla_NV,\lambda_kN\rangle\\
&=&\epsilon\lambda_k\phi\langle N,N\rangle=\lambda_k\phi
\end{eqnarray*}
and
$$\langle\overline\nabla_{e_k}N,\overline\nabla_{e_k}V\rangle=-\lambda_k\langle e_k,\overline\nabla_{e_k}V\rangle=-\lambda_k\phi,$$
so that
\begin{equation}\label{eq:IV}
\langle\overline\nabla_{e_k}\overline\nabla_NV,e_k\rangle+\langle N,\overline\nabla_{e_k}\overline\nabla_{e_k}V\rangle=0
\end{equation}
at $p$. On the other hand, since
$$[N,e_k](p)=(\overline\nabla_Ne_k)(p)-(\overline\nabla_{e_k}N)(p)=\lambda_ke_k(p),$$
it follows from (\ref{eq:IV}) that
\begin{eqnarray*}
\langle\overline R(N,e_k)V,e_k\rangle_p&=&\langle\overline\nabla_{e_k}\overline\nabla_NV-\overline\nabla_N\overline\nabla_{e_k}V+\overline\nabla_{[N,e_k]}V,e_k\rangle_p\\
&=&-\langle N,\overline\nabla_{e_k}\overline\nabla_{e_k}V\rangle_p-N\langle\overline\nabla_{e_k}V,e_k\rangle_p+\langle\overline\nabla_{\lambda_ke_k}V,e_k\rangle_p\\
&=&-\langle N,\overline\nabla_{e_k}\overline\nabla_{e_k}V\rangle_p-N(\phi)+\lambda_k\phi,
\end{eqnarray*}
and hence
\begin{equation}\label{eq:V}
\sum_k\langle N,\overline\nabla_{e_k}\overline\nabla_{e_k}V\rangle_p=-nN(\phi)+\epsilon\,nH\phi-\overline{\rm{Ric}}(N,V)_p
\end{equation}
Finally,
\begin{eqnarray*}
\overline{\rm{Ric}}(N,V)&=&\sum_l\alpha_l\overline{\rm{Ric}}(N,e_l)+\epsilon\,\eta\overline{\rm{Ric}}(N,N)\\
&=&\sum_{k,l}\alpha_l\langle\overline R(e_k,e_l)e_k,N\rangle+\epsilon\,\eta\overline{\rm{Ric}}(N,N),
\end{eqnarray*}
and
\begin{eqnarray*}
\langle\overline R(e_k,e_l)e_k,N\rangle_p&=&\langle\overline\nabla_{e_l}\overline\nabla_{e_k}e_k-\overline\nabla_{e_k}\overline\nabla_{e_l}e_k,N\rangle_p\\
&=&e_l\langle\overline\nabla_{e_k}e_k,N\rangle_p-\langle\overline\nabla_{e_k}e_k,\overline\nabla_{e_l}N\rangle_p-e_k\langle\overline\nabla_{e_l}e_k,N\rangle_p\\
&&+\langle\overline\nabla_{e_l}e_k,\overline\nabla_{e_k}N\rangle_p\\
&=&-e_l\langle e_k,\overline\nabla_{e_k}N\rangle_p+e_k\langle e_k,\overline\nabla_{e_l}N\rangle_p\\
&=&e_l(h_{kk})-e_k(h_{kl}),
\end{eqnarray*}
so that
$$\overline{\rm{Ric}}(N,V)_p=\sum_{k,l}\alpha_le_l(h_{kk})-\sum_{k,l}\alpha_le_k(h_{kl})+\epsilon\,\eta\overline{\rm{Ric}}(N,N)_p,$$
and it follows from (\ref{eq:V}) that
\begin{eqnarray}\label{eq:VI}
\sum_k\langle N,\overline\nabla_{e_k}\overline\nabla_{e_k}V\rangle_p&=&-nN(\phi)+\epsilon\,nH\phi-V^{\top}(\epsilon\,nH)\nonumber\\
&&+\sum_{k,l}\alpha_le_k(h_{kl})-\epsilon\,\eta\overline{\rm{Ric}}(N,N).
\end{eqnarray}

Substituting (\ref{eq:II}), (\ref{eq:III}) and (\ref{eq:VI}) into (\ref{eq:I}), one gets
the desired formula (\ref{eq:Laplacian formula_I}).

\end{proof}

\section{Semi-Riemannian warped products}\label{section:Semi-Riemannian warped products}

Let $M^n$ be a connected, $n$-dimensional oriented Riemannian manifold, $I\subset\mathbb R$ an interval and
$f:I\rightarrow\mathbb R$ a positive smooth function. In the product differentiable manifold
$\overline M^{n+1}=I\times M^n$, let $\pi_I$ and $\pi_M$ denote the projections onto the $I$ and $M$
factors, respectively.

A particular class of semi-Riemannian manifolds having conformal fields is the one obtained by
furnishing $\overline M$ with the metric
$$\langle v,w\rangle_p=\epsilon\langle(\pi_I)_*v,(\pi_I)_*w\rangle+f(p)^2\langle(\pi_M)_*v,(\pi_M)_*w\rangle,$$
where $\epsilon=-1$ or $\epsilon=1$ for all $p\in\overline M$ and all $v,w\in T_p\overline M$.
Indeed (cf.~\cite{Montiel1:99} and~\cite{Montiel2:99}),
the vector field
$$V=(f\circ\pi_I)\partial_t$$
is conformal and closed (in the sense that its dual $1-$form is closed), with conformal factor $\phi=f'$, where the
prime denotes differentiation with respect to $t\in I$. Such a space is called a semi-Riemannian {\em warped product},
and in what follows we shall write $\overline M^{n+1}=\epsilon I\times_fM^n$ to denote it.

If $\psi:\Sigma^n\rightarrow\overline M^{n+1}$ is a Riemannian immersion, with $\Sigma$ oriented by the unit
vector field $N$, one obviously has $\epsilon=\epsilon(\partial_t)=\epsilon(N)$.
The following result restates proposition~\ref{prop:Laplacian of conformal vector field} in this context,
in the spirit of~\cite{Alias2:06}.

\begin{proposition}\label{prop:Laplacian of conformal vector field for warped spaces}
Let $\overline M^{n+1}=\epsilon\,I\times_fM^n$. In the notations of
proposition~\ref{prop:Laplacian of conformal vector field}, if $\Sigma$ has constant mean curvature $H$, then
\begin{equation}\label{eq:Laplacian formula_II}
\Delta\eta=-\epsilon\,\eta\left\{{\rm Ric}(N^{\top},N^{\top})+(n-1)(\log f)''(1-\langle N,\partial_t\rangle^2)+|A|^2\right\}-\epsilon\,nHf'.
\end{equation}
where $\rm Ric$ denotes the Ricci tensor of $M$ and $N^{\top}=(\pi_M)_*N$.
\end{proposition}

\begin{proof} First of all,
$\eta=\langle V,N\rangle=f\langle N,\partial_t\rangle$, and
it thus follows from (\ref{eq:Laplacian formula_I}) that
$$\Delta\eta=-\epsilon\,\eta\left\{\overline{\rm Ric}(N,N)+|A|^2\right\}-n\left\{\epsilon Hf'+N(f')\right\}.$$
Now, $N(f')=\epsilon\,f''\langle N,\partial_t\rangle=\epsilon\,(f''/f)\eta$.
On the other hand, since $N=N^{\top}+\epsilon\langle N,\partial_t\rangle\partial_t$,
it follows from corollary $7.43$ of~\cite{O'Neill:83} that

\begin{eqnarray*}
\overline{\rm Ric}(N,N)&=&\overline{\rm Ric}(N^{\top},N^{\top})+\langle N,\partial_t\rangle^2\overline{\rm Ric}(\partial_t,\partial_t)\\
&=&{\rm Ric}(N^{\top},N^{\top})-\epsilon\langle N^{\top},N^{\top}\rangle\left\{\frac{f''}{f}+(n-1)\frac{(f')^2}{f^2}\right\}-\frac{nf''}{f}\langle N,\partial_t\rangle^2\\
&=&{\rm Ric}(N^{\top},N^{\top})-\left\{\frac{f''}{f}+(n-1)\frac{(f')^2}{f^2}\right\}-(n-1)\left(\frac{f'}{f}\right)'\langle N,\partial_t\rangle^2,
\end{eqnarray*}
where we used that $\langle N^{\top},N^{\top}\rangle=\epsilon(1-\langle N,\partial_t\rangle^2)$
in the last equality above.
\begin{eqnarray*}
\Delta\eta&=&-\epsilon\,\eta\left\{{\rm Ric}(N^{\top},N^{\top})-\left\{\frac{f''}{f}+(n-1)\frac{(f')^2}{f^2}\right\}-(n-1)\left(\frac{f'}{f}\right)'\langle N,\partial_t\rangle^2\right\}\\
&&-\epsilon\,\eta|A|^2-\epsilon\,n\left\{Hf'+\frac{f''}{f}\eta\right\}\\
&=&-\epsilon\,\eta\left\{{\rm Ric}(N^{\top},N^{\top})+(n-1)(\log f)''(1-\langle N,\partial_t\rangle^2)+|A|^2\right\}-\epsilon\,nHf'.
\end{eqnarray*}
\end{proof}

If $\psi:\Sigma^n\rightarrow\overline M^{n+1}$ is a Riemannian immersion as above, we
let $h=\pi_{I_{|\Sigma}}:\Sigma\rightarrow I$ denote the height function of $\Sigma$ with respect
to the unit vector field $\partial_t$. As far as we know, the following proposition appeared
for the first time in~\cite{Alias2:06}, as a special case of lemma 4.1; here we present
a direct proof of the particular case which is needed for the applications we have in mind.
We would like to thank professors L. Al\'{\i}as and A. G. Colares for having showed us the above mentioned
preprint, which enabled us to set this result.

\begin{proposition}\label{prop:Laplacian of h}
In the above notation,
\begin{equation}\label{eq:Laplacian of h}
\Delta h=(\log f)'(h)\{\epsilon\,n-|\nabla h|^2\}+\epsilon\,nH\langle N,\partial_t\rangle,
\end{equation}
where $H$ denotes the mean curvature of $\Sigma$ with respec to $N$.
\end{proposition}

\begin{proof}
Since $h=\pi_{I_{|\Sigma}}$, one has
\begin{eqnarray*}
\nabla h&=&\nabla(\pi_{I_{|\Sigma}})=(\overline\nabla\pi_I)^{\top}=\epsilon\partial_t^{\top}\\
&=&\epsilon\,\partial_t-\langle N,\partial_t\rangle N.
\end{eqnarray*}
where $\overline\nabla$ denotes the gradient with respect to the metric of the ambient space, and
$X^{\top}$ the tangential component of a vector field $X\in\mathcal X(\overline M)$ in $\Sigma$.
Now fix $p\in M$, $v\in T_pM$ and let $A$ denote the Weingarten map with respect to $N$. Write
$v=w+\epsilon\langle v,\partial_t\rangle\partial_t$, so that $w\in T_p\overline M$ is tangent to the
fiber of $\overline M$ passing through $p$. Therefore, by repeated use of the formulas of item (2) of proposition $7.35$
of~\cite{O'Neill:83}, we get
\begin{eqnarray*}
\overline\nabla_v\partial_t&=&\overline\nabla_w\partial_t+\epsilon\langle v,\partial_t\rangle\overline\nabla_{\partial_t}\partial_t=\overline\nabla_w\partial_t\\
&=&(\log f)'w=(\log f)'(v-\epsilon\langle v,\partial_t\rangle\partial_t),
\end{eqnarray*}
so that
\begin{eqnarray*}
\nabla_v\nabla h&=&\overline\nabla_v\nabla h-\epsilon\langle Av,\nabla h\rangle N\\
&=&\overline\nabla_v(\epsilon\,\partial_t-\langle N,\partial_t\rangle N)-\epsilon\langle Av,\nabla h\rangle N\\
&=&\epsilon(\log f)'w-v(\langle N,\partial_t\rangle)N+\langle N,\partial_t\rangle Av-\epsilon\langle Av,\nabla h\rangle N\\
&=&\epsilon(\log f)'w+(\langle Av,\partial_t\rangle-\langle N,\overline\nabla_v\partial_t\rangle)N+\langle N,\partial_t\rangle Av-\epsilon\langle Av,\nabla h\rangle N\\
&=&\epsilon(\log f)'w+(\langle Av,\partial_t^{\top}\rangle-\langle N,(\log f)'w\rangle)N+\langle N,\partial_t\rangle Av-\epsilon\langle Av,\nabla h\rangle N\\
&=&\epsilon(\log f)'w+\epsilon(\log f)'\langle v,\partial_t\rangle\langle N,\partial_t\rangle N+\langle N,\partial_t\rangle Av\\
&=&\epsilon(\log f)'\{v-\langle v,\partial_t\rangle(\epsilon\,\partial_t-\langle N,\partial_t\rangle N)\}+\langle N,\partial_t\rangle Av\\
&=&(\log f)'(\epsilon\,v-\epsilon\langle v,\partial_t^{\top}\rangle\nabla h)+\langle N,\partial_t\rangle Av\\
&=&(\log f)'(\epsilon\,v-\langle v,\nabla h\rangle\nabla h)+\langle N,\partial_t\rangle Av\\
\end{eqnarray*}
Now, fixing $p\in\Sigma$ and an orthonormal frame $\{e_i\}$ at $T_p\Sigma$, one gets
\begin{eqnarray*}
\Delta h&=&{\rm tr}(\nabla^2h)=\sum_{i=1}^n\langle\nabla_{e_i}\nabla h,e_i\rangle\\
&=&\sum_{i=1}^n\langle(\log f)'(\epsilon\,e_i-\langle e_i,\nabla h\rangle\nabla h)+\langle N,\partial_t\rangle Ae_i,e_i\rangle\\
&=&(\log f)'\{\epsilon\,n-|\nabla h|^2\}+\langle N,\partial_t\rangle\rm{tr}(A)\\
&=&(\log f)'\{\epsilon\,n-|\nabla h|^2\}+\epsilon\,nH\langle N,\partial_t\rangle.
\end{eqnarray*}
\end{proof}

Let us consider again a semi-Riemannian warped product $\overline M^{n+1}=\epsilon I\times_fM^n$.
For $t_0\in\mathbb R$, we orient the fiber $M_{t_0}^n=\{t_0\}\times M^n$ by using the unit normal vector
field $\partial_t$. According to proposition $1$ of~\cite{Montiel1:99} (see also proposition
$1$ of~\cite{Montiel2:99}), $M_{t_0}$ has constant mean curvature $-\epsilon f'(t_0)/f(t_0)$.
We are finally in position to define the objects of our main concern.

\begin{definition}\label{def:vertical graphs}
Let $\psi:\Sigma^n\rightarrow\overline M^{n+1}$ be a Riemannian immersion. We say that $\Sigma$ is a
{\em vertical graph} over the fiber $M_{t_0}^n$ if $\psi(x)=(u(x),x)$ for some smooth function
$u:M_{t_0}\rightarrow[0,+\infty)$.
\end{definition}

Three remarks are in order. First of all, if we let $h$ denote the height function associated to a vertical graph
over the fiber $M_{t_0}$, with corresponding function $u:M_{t_0}\rightarrow[0,+\infty)$, then one obviously has
$u=h\circ\psi-t_0$. Secondly, in the Lorentz case the condition that $\psi$ is Riemannian in the above definition
amounts to $|Du|<1$, where by $Du$ we mean the gradient of $u\circ\iota$ with respect to the metric of $M$,
where $\iota:M\rightarrow M_{t_0}$ is the canonical map (cf.~\cite{Montiel:03}, section $4$).
At last, our applications in the following sections all deal with semi-Riemannian warped products with warping
function $f(t)=e^t$. According to the dicussion preceeding the above definition, in this setting all fibers
have mean curvature $-\epsilon$,
and due to this fact we will assume that our vertical graphs are those over $M_0$, i.e., such that
$u=h\circ\psi\geq 0$. This agreement clarifies our exposition and does not imply in any loss of generality; indeed,
changing $u$ by $u+t_0$, all of the arguments to come can be easily adapted to vertical graphs over $M_{t_0}$.

\section{Vertical graphs in the Steady State space}\label{section:Vertical graphs in the steady state space}

In this section we consider a particular model of Lorentzian warped product, the
{\em Steady State space}, namely, the warped product
\begin{equation}\label{eq:steady state space}
\mathcal{H}^{n+1}=-\mathbb R\times_{e^t}\mathbb R^n.
\end{equation}
In Cosmology, this space corresponds to the steady state model of the
universe proposed by Bondi, Gold and Hoyle (cf.~\cite{Hawking:73}, p. $126$).

An alternative description of the Steady State space $\mathcal{H}^{n+1}$ (cf.~\cite{Montiel:03}) can be given
as follows. Let $\mathbb L^{n+2}$ denote the $(n+2)$-dimensional
Lorentz-Minkowski space ($n\geq 2$), that is, the real vector space $\mathbb R^{n+2}$,
endowed with the Lorentz metric
$$\left\langle v,w\right\rangle ={\displaystyle\sum\limits_{i=1}^{n+1}}v_{i}w_{i}-v_{n+2}w_{n+2},$$
for all $v,w\in\mathbb{R}^{n+2}$. We define the $\left(n+1\right)$-dimensional de Sitter
space $\mathbb S_1^{n+1}$ as the hyperquadric
$$\mathbb S_1^{n+1}=\left\{p\in L^{n+2};\left\langle p,p\right\rangle=1\right\}$$
of $\mathbb{L}^{n+2}$. From the above definition it is easy to show that the metric
induced from $\left\langle\,\,,\,\right\rangle $ turns $\mathbb S_1^{n+1}$ into a
Lorentz manifold with constant sectional curvature $1$. Moreover, for
$p\in\mathbb S_1^{n+1}$, we have
$$T_p\mathbb S_1^{n+1}=\left\{v\in\mathbb L^{n+2};\left\langle v,p\right\rangle=0\right\}.$$

Let $a\in\mathbb L^{n+2}$ be a nonzero null vector of the null cone with vertex in the origin, such that
$\left\langle a,e_{n+2}\right\rangle >0$, where $e_{n+2}=(0,\ldots,0,1)$.
It can be shown that the open region
$$\left\{p\in\mathbb S_1^{n+1};\left\langle p,a\right\rangle>0\right\}$$
of the de Sitter space $\mathbb S_1^{n+1}$ is isometric to $\mathcal H^{n+1}$.
Therefore, as a subset of $\mathbb S_1^{n+1}$, the boundary of $\mathcal H^{n+1}$
is the null hypersurface
$$\left\{p\in\mathbb S_1^{n+1};\left\langle p,a\right\rangle =0\right\}.$$

Back to the warped product model of $\mathcal H^{n+1}$, if $\psi:\Sigma^n\rightarrow\mathcal H^{n+1}$
is a spacelike hypersurface oriented by the timelike unit vector field $N$ such that
$\langle N,\partial_t\rangle<0$, the {\em hyperbolic angle} $\theta$ of $\psi$ is
the smooth function $\theta:\psi(\Sigma)\rightarrow[0,+\infty)$ such that
\begin{equation}\label{eq:hyperbolic angle}
\cosh\theta=-\langle N,\partial_t\rangle\geq 1.
\end{equation}

In the following result, the right hand side of (\ref{eq:Lorentz growth of h}) must be interpreted as $+\infty$ when
$\cosh\theta=1$.

\begin{theorem}\label{thm:de Sitter application}
Let $\psi:\Sigma^n\rightarrow\mathcal H^{n+1}$ be a complete spacelike vertical graph in the $(n+1)$-dimensional Steady
State space, with constant mean curvature $H\geq 1$. If
\begin{equation}\label{eq:Lorentz growth of h}
h\leq-\log(\cosh\theta-1),
\end{equation}
then:
\begin{enumerate}
\item[$(a)$] $H=1$ on $\Sigma$.
\item[$(b)$] The scalar curvature $R$ of $\Sigma$ is nonnegative and not globally bounded away from zero.
\end{enumerate}
\end{theorem}

\begin{proof}
Let $g:\Sigma\rightarrow\mathbb R$ be defined by $g=-e^h-\eta$. It follows easily from (\ref{eq:hyperbolic angle}) and
the definition of $h$ that $g\geq 0$ on $\Sigma$. On the other hand, our hypothesis on the growth of $h$ assures that
$g\leq 1$ on $\Sigma$.

A straightforward computation gives us $\Delta e^h=e^h\{|\nabla h|^2+\Delta h\}$. Moreover, since the Riemannian
fiber of $\mathcal H^{n+1}$ is $\mathbb R^n$, by computing the Laplacian of $g$
with the aid of propositions~\ref{prop:Laplacian of conformal vector field for warped spaces}
and~\ref{prop:Laplacian of h} we get
\begin{eqnarray*}
\Delta g&=&-\Delta e^h-\Delta\eta\\
&=&-e^h\{|\nabla h|^2+\Delta h\}-\Delta\eta\\
&=&ne^h\{1+H\langle N,\partial_t\rangle\}-\eta|A|^2-nHe^h.
\end{eqnarray*}

Now, let $S_2$ denote the second elementary symmetric function on the eigenvalues of $A$, and $H_2=2S_2/n(n-1)$ denote
the mean value of $S_2$. Elementary algebra gives
$$|A|^2=n^2H^2-n(n-1)H_2,$$
which put into the above formula gives, after a little more algebra,
\begin{eqnarray}\label{eq:key expression for the Laplacian of g in de Sitter}
\Delta g&=&n(H-1)\{-e^h-H\eta\}-n(n-1)(H^2-H_2)\eta\nonumber\\
&\geq&n(H-1)g+n(n-1)(H^2-H_2),
\end{eqnarray}
where for the inequality we used that $-\eta\geq e^h\geq 1$.\\

\noindent $(a)$ Suppose, by contradiction, that $H>1$. Since
$0\leq g\leq 1$ and (from the Cauchy-Schwarz inequality) $H^2-H_2\geq 0$, we get
$$\Delta g\geq n(H-1)g^2.$$
Now let ${\rm Ric}_{\Sigma}$ denote the Ricci curvature of $\Sigma$; by applying Gauss' equation, we get the estimate
\begin{equation}\label{eq:estimate of the Ricci curvature}
{\rm Ric}_{\Sigma}\geq(n-1)-\frac{n^2H^2}{4},
\end{equation}
so that we are in position to apply lemma $3$ of~\cite{Akutagawa:87} to conclude that $g\equiv 0$. Thus,
$\eta\equiv-e^h$, so that $\langle N,\partial_t\rangle\equiv-1$, i.e., $\psi(\Sigma)$ is a slice of $\mathcal H$.
However, such a slice has constant mean curvature $1$, and we arrive at a contradicion. Thus $H=1$.\\

\noindent $(b)$ Back to (\ref{eq:key expression for the Laplacian of g in de Sitter}), we obtain
$$\Delta g\geq n(n-1)(1-H_2)=R\geq 0,$$
where we used Gauss' equation once more to get the last equality, and $H^2-H_2\geq 0$ to get the sign for $R$.

Hence, if there exists $\alpha>0$ such
that $R\geq\alpha$ on $\Sigma$, from the above we could derive the inequality
$$\Delta g\geq\alpha g^2,$$
which once more would give us $g\equiv 0$, so that $\psi(\Sigma)$ would also be a slice. However, such a slice is
isometric to $\mathbb R^n$, thus having scalar curvature $R\equiv 0$. We, therefore, have got another contradiction.
\end{proof}

\begin{remark}
It is easy to see that hypothesis (\ref{eq:Lorentz growth of h}) on the growth of $h$ is
implied by the estimate
$$|\nabla h|\leq e^{-h/2}$$
for the gradient of $h$, which in turn is taken as a natural one in the literature
(see, for instance, corollary $16.6$ of~\cite{Gilbarg:83}).
\end{remark}

\begin{remark}
As a consequence of Bonnet-Myers theorem, a complete spacelike
hypersurface $\psi:\Sigma^n\rightarrow\mathcal H^{n+1}$ having (not necessarily constant)
mean curvature $H$, such that $|H|\leq\varrho<2\sqrt{n-1}/n$ ($\varrho$ constant), has to
be compact; in fact, for such a bound on $H$, equation
(\ref{eq:estimate of the Ricci curvature}) would give us
$${\rm Ric}_{\Sigma}\geq(n-1)-n^2\rho^2/4>0.$$
However, since $\psi(\Sigma)$ is a graph over $\mathbb R^n$, it cannot be compact. Therefore, since
$2\sqrt{n-1}/n\leq 1$ for $n\geq 2$, in a certain sense it is natural to restrict attention to $H\geq 1$.
\end{remark}

As a consequence of the previous result, we have the following Bernstein-type theorem in $\mathcal H^3$:

\begin{theorem}\label{thm:from Akutagawa}
Let $\psi:\Sigma^2\rightarrow\mathcal H^3$ be a complete spacelike vertical graph in the $3$-dimensional Steady
State space, with constant mean curvature $H\geq 1$. If
$$h\leq-\log(\cosh\theta-1),$$
then $\psi(\Sigma)$ is a slice of $\mathcal H^3$.
\end{theorem}

\begin{proof}
From the previous result, $H=1$ on $\Sigma$. Now apply the main theorem of~\cite{Akutagawa:87} and the
classification of umbilical hypersurfaces of the de Sitter space, cf. example $1$ of~\cite{Montiel:88}.
\end{proof}

We can also apply the result of Proposition~\ref{prop:Laplacian of h} to prove yet another Bernstein-type theorem
for complete surfaces (not necessarily graphs) of the $3-$dimensional Steady State space.

\begin{theorem}\label{thm:de Sitter Bernstein-type theorem}
Let $\psi:\Sigma^2\rightarrow\mathcal H^3$ be a Riemannian immersion of a complete surface
of nonnegative Gaussian curvature $K_{\Sigma}$, with constant mean curvature $H\geq 1$. If
\begin{equation}\label{eq:gradient bound in De Sitter}
|\nabla h|^2\leq H^2-1,
\end{equation}
then $\psi\left(\Sigma\right)$ is a slice of $\mathcal H^3$.
\end{theorem}

\begin{proof}
By applying the result of Proposition~\ref{prop:Laplacian of h}, we get
\begin{eqnarray*}
\triangle e^{-h} &=& e^{-h}\left\{|\nabla h|^2-\triangle h\right\} \\
&=& 2e^{-h}\left\{|\nabla h|^2+1+H\langle N,\partial_t\rangle \right\}.
\end{eqnarray*}

On the other hand, since $|\nabla h|^2=\langle N,\partial_t\rangle^2-1$,
hypothesis (\ref{eq:gradient bound in De Sitter}) is equivalent to
$$|\nabla h|^2+1+H\langle N,\partial_t\rangle\leq 0,$$
so that the function $e^{-h}$ is a superharmonic positive function on $\Sigma$.
However, a classical result due to Huber~\cite{Huber:57} assures that complete surfaces of non-negative
Gaussian curvature must be parabolic; therefore, $h$ is constant on $\Sigma$, i.e., $\psi\left(\Sigma\right)$
is a slice.
\end{proof}

%
%
%
%

\section{Vertical graphs in the Hyperbolic space}\label{section:Vertical graphs in the Hyperbolic space}

In this section, instead of the more commonly used half-space model for the
$(n+1)-$dimensional Hyperbolic space, we consider the warped product model
$$\mathbb H^{n+1}=\mathbb R\times_{e^t}\mathbb R^n.$$
An explicit isometry between these two models can be found at~\cite{Alias:06}, from
where it can easily be seen that the fibers $M_{t_0}=\{t_0\}\times\mathbb R^n$ of the warped product
model are precisely the horospheres of $\mathbb H^{n+1}$. Moreover, according to the last paragraph of
section~\ref{section:Semi-Riemannian warped products}, these have constant mean curvature $1$ if we
take the orientation given by the unit normal vector field $N=-\partial_t$.

Another useful model for $\mathbb H^{n+1}$ is (following the notation of the previous section) the so-called
{\em Lorentz model}, obtained by furnishing the hyperquadric
$$\{p\in\mathbb L^{n+2};\,\langle p,p\rangle=-1,\,p_{n+2}>0\}$$
with the (Riemannian) metric induced by the Lorentz metric of $\mathbb L^{n+2}$.
In this setting, if $a\in\mathbb L^{n+2}$ denotes a fixed null vector as in the beginning of the previous
section, a typical horosphere is
$$L_{\tau}=\{p\in\mathbb H^{n+1};\,\langle p,a\rangle=\tau\},$$
where $\tau$ is a positive real number. A straightforward computation shows that
$$\xi_p=p+\frac{1}{\tau}a\in\mathcal H^{n+1}$$
is a unit normal vector field along $L_{\tau}$, with respect to which $L_{\tau}$ has mean curvature $-1$
(cf.~\cite{Lopez:99}). Therefore, any isometry $\Phi$ between the warped product and Lorentz models
of $\mathbb H^{n+1}$ must carry $(\partial_t)_q$ to $\Phi_*(\partial_t)=\xi_{\Phi(q)}$.

If $\psi:\Sigma^n\rightarrow\mathbb H^{n+1}$ is a vertical graph over $\mathbb R^n$, we orient $\Sigma$ by
choosing a unit normal vector field $N$ such that $\eta=\langle N,V\rangle<0$, and hence $-e^h\leq\eta<0$.
Following the discussion of the previous paragraph, it is natural to consider the {\em Lorentz Gauss map} of $\Sigma$
with respect to $N$ as given by
$$\begin{array}{ccc}
\Sigma^n&\rightarrow&\mathcal H^{n+1}\\
p&\mapsto&-\Phi_*(N_p)
\end{array}$$

We are finally in position to state and prove, in the Hyperbolic setting,
analogues of two of the results of the previous section, starting with theorem~\ref{thm:de Sitter application}.

\begin{theorem}\label{thm:Hyperbolic application}
Let $\Sigma$ be a complete Riemannian manifold with Ricci curvature globally bounded from below, and
$\psi:\Sigma^n\rightarrow\mathbb H^{n+1}$ be a vertical graph in the $(n+1)$-dimensional hyperbolic
space, with constant mean curvature $0\leq H\leq 1$. If
\begin{equation}\label{eq:growth of h}
h\leq-\log(1+\langle N,\partial_t\rangle),
\end{equation}
then:
\begin{enumerate}
\item[$(a)$] $H=1$ on $\Sigma$.
\item[$(b)$] If the closure of the image of the Lorentz Gauss map of $\psi$ with respect to $N$
is contained in $\mathcal H^{n+1}$, then the scalar curvature $R$ of $\Sigma$ is nonpositive and not
globally bounded away from zero.
\end{enumerate}
\end{theorem}

\begin{proof}
Let $g:\Sigma\rightarrow\mathbb R$ be defined by $g=e^h+\eta$. The definition of $h$, together with Cauchy-Schwarz
inequality, gives us $g\geq 0$ on $\Sigma$; on the other hand, our hypothesis on the growth of $h$ assures that
$g\leq 1$ on $\Sigma$.

A straightforward computation gives us $\Delta e^h=e^h\{|\nabla h|^2+\Delta h\}$. Moreover, since the Riemannian
fiber of $\mathbb H^{n+1}$ is $\mathbb R^n$, by computing the Laplacian of $g$
with the aid of propositions~\ref{prop:Laplacian of conformal vector field for warped spaces}
and~\ref{prop:Laplacian of h} we get
\begin{eqnarray*}
\Delta g&=&\Delta e^h+\Delta\eta\\
&=&e^h\{|\nabla h|^2+\Delta h\}+\Delta\eta\\
&=&ne^h\{1+H\langle N,\partial_t\rangle\}-\eta|A|^2-nHe^h.
\end{eqnarray*}

Now, let $S_2$ denote the second elementary symmetric function on the eigenvalues of $A$, and $H_2=2S_2/n(n-1)$ denote
the mean value of $S_2$. Elementary algebra gives
$$|A|^2=n^2H^2-n(n-1)H_2,$$
which put into the above formula gives, after a little more algebra,
\begin{eqnarray}\label{eq:key expression for the Laplacian of g}
\Delta g&=&n(1-H)\{e^h+H\eta\}-n(n-1)(H^2-H_2)\eta\nonumber\\
&=&n(1-H)g-n(n-1)(H^2-H_2)\eta.
\end{eqnarray}

\noindent $(a)$
Suppose, by the sake of contradiction, that $H<1$ on $\Sigma$. Since $0\leq g\leq 1$, $-\eta>0$ and $H^2-H_2\geq 0$
(from Cauchy-Schwarz inequality), we get
$$\Delta g\geq n(1-H)g^2.$$
Thus, from our hiypothesis on the Ricci curvature of $\Sigma$ we are in position to apply lemma $3$
of~\cite{Akutagawa:87} to conclude that $g\equiv 0$, which is the same as
$\langle N,\partial_t\rangle\equiv-1$. Therefore, $\psi(\Sigma)$ is a horosphere of $\mathbb H^{n+1}$.
However, such a horosphere has constant mean curvature $1$, and we reached a contradicion.\\

\noindent $(b)$ Back to (\ref{eq:key expression for the Laplacian of g}), we get
$$\Delta g=n(n-1)(H_2-1)\eta=R\eta\geq R\langle N,\partial_t\rangle,$$
where Gauss' equation was applied for the last equality and we used that $\eta<0$ and $H_2-H^2\leq 0$ for the
last inequality. The condition on the Lorentz Gauss map of $\Sigma$ amounts to the existence of a real
number $\beta>0$ such that $\langle-N,\partial_t\rangle\geq\beta$ on $\Sigma$. Therefore, if there existed a positive
real number $\alpha$ such that $R\leq-\alpha$ on $\Sigma$, we would get from $0\leq g\leq 1$ that
$$\Delta g\geq-R\langle-N,\partial_t\rangle\geq\alpha\beta g^2,$$
so that applying lemma $3$ of~\cite{Akutagawa:87} once more would give us $g\equiv 0$. However, horospheres of
$\mathbb H^{n+1}$ are isometric to $\mathbb R^n$, thus having scalar curvature identically $0$, which is a
contradiction.
\end{proof}

We close this paper with an analogue of theorem~\ref{thm:de Sitter Bernstein-type theorem}
for the Hyperbolic space.

\begin{theorem}\label{thm:Hyperbolic Bernstein-type theorem}
Let $\psi:\Sigma^2\rightarrow\mathbb H^3$ be a complete vertical graph with nonnegative
Gaussian curvature $K_{\Sigma}$ and constant mean curvature $\frac{\sqrt{2}}{2}\leq H\leq1$.
If
\begin{equation}\label{eq:Hyperbolic gradient estimate}
|\nabla h|^2\leq 1-H^2,
\end{equation}
then $\psi\left(\Sigma\right)$ is a horosphere of $\mathbb H^3$.
\end{theorem}

\begin{proof}
Once more from proposition~\ref{prop:Laplacian of h}, we have
$$\triangle e^{-h}=2e^{-h}\left\{|\nabla h|^2-1-H\langle N,\partial_t\rangle \right\}.$$

On the other hand, since $|\nabla h|^2=1-\langle N,\partial_t\rangle^2$ and
$\langle N,\partial_t\rangle$ does not change sign, hypothesis
(\ref{eq:Hyperbolic gradient estimate}) is equivalent to
$$|\nabla h|^2-1-H\langle N,\partial_t\rangle\leq 0,$$
so that $e^{-h}$ is a superharmonic and positive on $\Sigma^2$. Hence, as in the proof of
theorem~\ref{thm:de Sitter Bernstein-type theorem}, $h$ is constant on $\Sigma^2$, i.e.,
$\psi(\Sigma)$ is a horosphere.
\end{proof}

\begin{remark}
Since Gauss's equation gives
$$K_{\Sigma}=2H^2-1-\frac{1}{2}|A|^2,$$
the assumption $K_{\Sigma}\geq 0$ forces one to restric attention to the case
$H\geq\frac{\sqrt 2}{2}$.
\end{remark}

\begin{remark}
Under the assumption of properness for $\psi$, a result similar to
the above can be found in~\cite{Alias:06}.
\end{remark}

\end{document}